\documentclass{article}
\usepackage{graphicx}
\usepackage{amsmath}
\usepackage{amsfonts}
\usepackage{amssymb}
\newtheorem{theorem}{Theorem}

\newtheorem{condition}{Condition}

\newtheorem{corollary}{Corollary}

\newtheorem{definition}{Definition}
\newtheorem{example}{Example}

\newtheorem{lemma}{Lemma}

\newtheorem{proposition}{Proposition}
\newtheorem{remark}{Remark}

\newenvironment{proof}[1][Proof]{\textbf{#1.} }{\ \rule{0.5em}{0.5em}}

\begin{document}

\title{On weight graphs for nilpotent Lie algebras I }
\author{Jos\'{e} Mar\'{\i}a Ancochea Berm\'{u}dez\thanks{%
research supported by the D.G.Y.C.I.T project PB98-0758} \and Otto Rutwig
Campoamor Stursberg\thanks{corresponding author : e-mail : rutwig@nfssrv.mat.ucm.es} \\
Departamento de Geometr\'{\i}a y Topolog\'{\i}a\\
Fac. CC. Matem\'{a}ticas Univ. Complutense\\
28040 Madrid ( Spain )}
\date{}
\maketitle
\begin{abstract}
We introduce the concept of weight graph for the weight system $P\frak{g}\left(T\right)$ of a finite dimensional nilpotent Lie algebra $\frak{g}$ and analyze the necessary conditions for a $\left(p,q\right)$-graph to be a weight graph for some $\frak{g}$.\newline AMS Subj. Clas. 05 C90, 17B30\newline \textit{keywords : Lie algebra, nilpotent, weight, graph}    
\end{abstract}

\section{Introduction}
Graph theory is probably the most "rediscovered" theory in mathematics. To the pure graph theory a lot of mathematical and non-mathematical disciplines have developed a graph theory more than once. From the study of statistical mechanics [Uh] to problems of algebraic topology, graph theory has shown its wideness of applicability. This paper is the first from a series dedicated to the developement of graph theory applied to the study of weight systems for finite dimensional nilpotent Lie algebras. There are in fact relations between combinatorics and Lie algebras, such as the statement of the four color problem or the study of classical affine Lie algebras [Pr]. However, the developement we propose here is quite different. We pretend to establish a graph theory that allows to study the torus of derivations for certain nilpotent Lie algebras, and which, conversely, help to decide wheter a given graph can appear as the weight graph of some nilpotent Lie algebra. For elementary properties and definitions about Lie algebras, we refer to reference [Hu], while for questions about graph theory we refer to [Ha]. Here we shall only consider finite undirected graphs without loops or multiple edges.

\subsection{Preliminaries and notations}

Let $\frak{g}$ be a finite dimensional complex nilpotent Lie algebra. Let
$Der\left(  \frak{g}\right)  $ be its Lie algebra of derivations. A torus $T$
over $\frak{g}$ is an abelian subalgebra of $Der\left(  \frak{g}\right)  $
consisting of semi-simple endomorphisms. Clearly the torus $T$ induces a
natural representation [Wi] on the Lie algebra $\frak{g}$, such that this
decomposes as
\[
\frak{g}=\bigsqcup_{\alpha\in T^{\ast}}\frak{g}_{\alpha}%
\]
where $T^{\ast}=Hom_{\mathbb{C}}\left(  T,\mathbb{C}\right)  $ and
$\frak{g}_{\alpha}=\left\{  X\in\frak{g\;}|\;\left[  t,X\right]
=\alpha\left(  t\right)  X\;\;\forall t\in T\right\}  $ is the weight space
corresponding to the weight $\alpha$. If the torus is maximal for the
inclusion relation, as the tori are conjugated [Mo], its common dimension is a
numerical invariant of $\frak{g}$ called the rank and denoted by $r\left(
\frak{g}\right)  $. Following Favre [Fa], we call
\[
R\frak{g}\left(  T\right)  =\left\{  \alpha\in T^{\ast}\;|\;\frak{g}_{\alpha
}\neq0\right\}
\]
the set of weights for the representation of $T$ over $\frak{g}$ and
\[
P\frak{g}\left(  T\right)  =\left\{  \left(  \alpha,d\alpha\right)
\;|\;\alpha\in R\frak{g}\left(  T\right)  ,\;d\alpha=\dim\frak{g}_{\alpha
}\right\}
\]

\begin{definition}
Let $\frak{g}$ be a nilpotent Lie algebra and  $T$ a maximal torus of
derivations. Then $R\frak{g}\left(  T\right)  $ is called a weight system for
$\frak{g}$.
\end{definition}

\begin{remark}
Clearly, given two maximal tori $T_{1}$ and $T_{2}$ over $\frak{g}$ such that
the induced representation over $\frak{g}$ is equivalent, it results that the
corresponding weight systems $P\frak{g}\left(  T_{1}\right)  $ and
$P\frak{g}\left(  T_{2}\right)  $ are also equivalent. In fact, the
equivalence class of a weight system constitutes an invariant of the algebra
[Fa]. We are mainly interested on the weight sets and the relations among the
individual weights. Though the weights can be chosen quite arbitrarily, the
relations they satisfy are preserved.
\end{remark}

\begin{condition}
Given the $n$-dimensional Lie algebra $\frak{g}$, we have $r\frak{g}\left(
T\right)  =\left\{  \alpha_{1},..,\alpha_{r}\right\}  $, where $r\leq n$. Here
we will only consider the case of $r=n$, i. e., the algebra decomposes through
the action of the torus into one dimensional weight subspaces. Observe that in
this case we can identify both the set of weights and the weight system, as
the multiplicities of weigths are one. We impose an additional condition,
namely, if there exist weights $\alpha_{i},\alpha_{j}\in R\frak{g}\left(
T\right)  $ such that $\alpha_{i}+\alpha_{j}\in R\frak{g}\left(  T\right)  $,
then there are vectors $X_{i}\in\frak{g}_{\alpha_{i}},X_{j}\in\frak{g}%
_{\alpha_{j}}$ such that $\left[  X_{i},X_{j}\right]  \neq0$. From now on,
whenever we speak in this work about a weight system, we will suppose that it
satisfies the preceding conditions. Although this condition seems to be very specific, it is not in fact too restrictive, as the nilpotent Lie algebras of maximal nilpotence index satisfy this condition [Ve]. These algebras are a very important class and have constituted the principal study object in nilpotent Lie algebras over the last thirty years.
\end{condition}

\begin{remark}
It is well known that the weight system alone does not determine the Lie
algebra law. There exist parametrized families with the same underlying
system. The preceding condition excludes exactly those conditions for which
the structure constants depending on the parameters vanish. It follows that
only a set of solutions of measure zero has been excluded.
\end{remark}

Suppose that $R\frak{g}\left(  T\right)  =\left\{  \alpha_{1},..,\alpha
_{n}\right\}  $ is the weight system of $\frak{g}$. To this weight system we
can associate the follwoing diagramm : let $V\left(  G\right)  =\left\{
p_{1},..,p_{n}\right\}  $ be the points, where $p_{i}$ corresponds to the weights $\alpha_{i}$, for all $i$.
We say that $p_{i}$ joins $p_{j}$ if $\alpha_{i}+\alpha_{j}\in
R\frak{g}\left(  T\right)  $. The corresponding graph is denoted by $G\left(
R\frak{g}\left(  T\right)  \right)  $.

\begin{lemma}
Let $R\frak{g}\left(  T\right)  $ be a weight system for the nilpotent Lie
algebra $\frak{g}$. Then $G\left(  R\frak{g}\left(  T\right)  \right)  $
contains at least an isolated point.
\end{lemma}

\begin{proof}
As the algebra is finite dimensional, it follows from the weight space
decomposition that there exists at least one weight $\gamma\in R\frak{g}%
\left(  T\right)  $ which is ''extreme'' in the following sense :%
\[
\alpha+\gamma\notin R\frak{g}\left(  T\right)  ,\;\forall \alpha\in R\frak{g}%
\left(  T\right)
\]
Then for any vector $X\in\frak{g}_{\gamma}$ we have $ad_{\frak{g}}X\equiv0$,
where $ad_{\frak{g}}$ denotes the adjoint operator in $\frak{g}$. Thus the
vector is central, and the corresponding vertex of $G\left(  R\frak{g}\left(
T\right)  \right)  $ is isolated.
\end{proof}

\begin{corollary}
The graph $G\left(  R\frak{g}\left(  T\right)  \right)  $ contains $\left(
\dim Z\left(  \frak{g}\right)  \right)  $ isolated points, where $Z\left(
\frak{g}\right)  $ is the center of $\frak{g}$.
\end{corollary}

\section{The weight graph $G=\overline{G\left(  R\frak{g}%
\left(  T\right)  \right)  }$}

According to [Ha], we denote by $\overline{G\left(  R\frak{g}\left(  T\right)
\right)  \text{ }}$ the complementary graph to $G\left(  R\frak{g}\left(
T\right)  \right)  $. It follows immediately from the lemma that
$\overline{G\left(  R\frak{g}\left(  T\right)  \right)  }$ is a connected graph.

\begin{definition}
Let $R\frak{g}\left(  T\right)  $ be a weight system of $\frak{g}$. Then the
graph $\overline{G\left(  R\frak{g}\left(  T\right)  \right)  }$ is called the
weight graph of $\frak{g}$.
\end{definition}

As the weight graphs are connected, we can define a metric on them. As usual,
we define the distance of two points $p_{i}$ and $p_{j}$ of $V\left(
G\right)  $ as the lenght of a shortest path joining them, noted $d\left(
p_{i},p_{j}\right)  $.

\begin{proposition}
Let $G$ be the weight graph for the nilpotent Lie algebra $\frak{g}$. Then,
for any two points $p_{i}$ and $p_{j}$ of $V\left(  G\right)  $, we have
\[
d\left(  p_{i},p_{j}\right)  \leq2
\]
\end{proposition}

\begin{proof}
As $G=\overline{G\left(  R\frak{g}\left(  T\right)  \right)  }$ for a weight
system $R\frak{g}\left(  T\right)  $ of $\frak{g}$, we know that $\overline
{G}$ contains at least one isolated point, call it $p_{0}$. Thus, in $G,$
$p_{0}$ is adjacent to any other point. If $p_{0}\neq p_{i}, p_{0}\neq P_{j}$ for $p_{i},p_{j}\in
V\left(  G\right)  $, then the path $p_{i}p_{0}p_{j}$ is a geodesic joining them. Thus
$d\left(  p_{i},p_{j}\right)  =2$. If one of the \ points equals $p_{0}$, then the
distance is one.
\end{proof}

Recall that for a given graph $G$, the n$^{th}$ power of $G$, $G^{n}$, is
defined as the graph whose points are those of $G$, and where $u,v$ are
adjacent if $d\left(  u,v\right)  \leq n$.

\begin{theorem}
Let $G$ be the weight graph for the nilpotent Lie algebra $\frak{g}$. Then
$G=G^{n}$ for any $n\geq2$.
\end{theorem}

This result gives a first strong reduction for a graph being the weight graph
of a nilpotent Lie algebra. Suppose that $G=\overline{G\left(  R\frak{g}%
\left(  T\right)  \right)  }$ is a $\left(  p,q\right)  $-graph, i.e., $G$ has
$p$ points and $q$ lines. This last number is lower or equal than $\left(
\begin{tabular}
[c]{l}%
$p$\\
$2$%
\end{tabular}
\right)  $. For $\overline{G}=G\left(  R\frak{g}\left(  T\right)  \right)  $
we have an isolated point, thus its number of lines is
\[
t\leq\left(
\begin{array}
[c]{c}%
p\\
2
\end{array}
\right)  -\left(  p-1\right)
\]
and it follows

\begin{proposition}
Let $G=\overline{G\left(  R\frak{g}%
\left(  T\right)  \right)  }$ be the weight graph for a nilpotent Lie algebra $\frak{g}$. Then the
number $q$ of lines in $G$ satisfies
\[
q\geq\left(  p-1\right)
\]
\end{proposition}

\begin{remark}
This lower bound for the number of lines of a weight graph is not very good,
as it does not consider the weight character of its points. The next step is
to obtain a more appropiate bound.
\end{remark}

\begin{lemma}
For a $\left(  p,q\right)$ graph $G=G\left(  R\frak{g}%
\left(  T\right)  \right)$ we have
\[
q\leq\sum_{j=1}^{\left[  \frac{p}{2}\right]  }\left(  p-2j\right)
\]
where $\left[  \frac{p}{2}\right]  $ denotes the integer part of $\frac{p}{2}$.
\end{lemma}

\begin{proof}
As we have imposed $\dim\frak{g}_{\alpha}=1$ for any weight $\alpha\in
R\frak{g}\left(  T\right)  $, we can reorder the weights of $\frak{g}$ in such
manner that a relations $\alpha_{i}+\alpha_{j}=\alpha_{k}$ corresponds to a
sum $i+j=k$, where $k\leq p$. Thus the maximal number of sums of weights
equals the number of possibilities $i+j=k$ with $k\leq p$ and $1\leq i<j$. It
is easily seen that this number is precisely $\sum_{j=1}^{\left[  \frac{p}%
{2}\right]  }\left(  p-2j\right)  $.
\end{proof}

\begin{remark}
Observe that this bound has been obtained indepently of graph theory. It is
only based in the principal property of weights, namely, that certain sums of
them give another weight.
\end{remark}

For $p\geq3$ let us consider the function $f\left(  p\right)  =\left(
\begin{array}
[c]{c}%
p\\
2
\end{array}
\right)  -\left(  p-1\right)  -%
{\displaystyle\sum_{j=1}^{\left[  \frac{p}{2}\right]  }}
\left(  p-2j\right)  $. This function measures how far is the first estimation
of lines of $G$ when compared with the one of the lemma. For $p=3$ we have
$f\left(  3\right)  =0$, while for $p=4$ we have $f\left(  4\right)  =1$.

\begin{lemma}
For $p\geq4$ we have $f\left(  p\right)  >0$.
\end{lemma}

\begin{proof}
The assertion is obviously true for $p=4$. Suppose it is also for $p>4$. Then
we have
\begin{align*}
f\left(  p+1\right)   &  =\left(
\begin{array}
[c]{c}%
p+1\\
2
\end{array}
\right)  -\left(  p+1\right)  -1-\sum_{j=1}^{\left[  \frac{p+1}{2}\right]
}\left(  p+1-2j\right) \\
&  =\left(
\begin{array}
[c]{c}%
p\\
2
\end{array}
\right)  -\sum_{j=1}^{\left[  \frac{p+1}{2}\right]  }\left(  p+1-2j\right)
\end{align*}

\begin{enumerate}
\item  if $p\equiv1\;\left(  \operatorname{mod}\;2\right)  $ then we have
$\left[  \frac{p+1}{2}\right]  =\left[  \frac{p}{2}\right]  +1$, and as
\[
\left(  p+1\right)  -2\left(  \left[  \frac{p}{2}\right]  +1\right)  =0,
\]
we have
\begin{align*}
f\left(  p+1\right)   &  =\left(
\begin{array}
[c]{c}%
p\\
2
\end{array}
\right)  -\sum_{j=1}^{\left[  \frac{p}{2}\right]  }\left(  p-2j\right)
-\left[  \frac{p}{2}\right] \\
&  =f\left(  p\right)  +\left(  p-1-\left[  \frac{p}{2}\right]  \right)
\end{align*}
As $p>4$ and $f\left(  p\right)  >0$ , it follows $f\left(  p+1\right)  >0.$

\item  if $p\equiv0\;\left(  \operatorname{mod}\;2\right)  $ then $\left[
\frac{p+1}{2}\right]  =\left[  \frac{p}{2}\right]  $. In this case
\[
f\left(  p+1\right)  =\left(
\begin{array}
[c]{c}%
p\\
2
\end{array}
\right)  -\sum_{j=1}^{\left[  \frac{p}{2}\right]  }\left(  p+1-2j\right)
-\left[  \frac{p}{2}\right]
\]
and as before $f\left(  p+1\right)  >0$.
\end{enumerate}
\end{proof}

In fact we have proven more than necessary, namely that for $p\geq3$ we have
\[
f\left(  p+1\right)  -f\left(  p\right)  =p-1-\left[  \frac{p}{2}\right]
\]

\begin{proposition}
Let $p\geq3$ and $k\geq1$. Then
\[
f\left(  p+k\right)  -f\left(  p\right)  =k\left(  p-1\right)  +\frac{k\left(
k-1\right)  }{2}-\sum_{j=0}^{k-1}\left[  \frac{p+j}{2}\right]
\]
\end{proposition}

This property is easily proven with induction over $k$. The next result shows
how rapidly the differences of the two bounds given increase :\newpage

\begin{proposition}
Let $p\geq8$. Then $f\left(  p\right)  \geq p+1$
\end{proposition}

\begin{proof}
Let $p=8+k$ with $k\geq1$. From the previous result we know that
\begin{align*}
f\left(  8+k\right)  -f\left(  8\right)   &  =7k+\frac{k\left(  k-1\right)
}{2}-\sum_{j=0}^{k-1}\left[  \frac{8+j}{2}\right] \\
&  =9+k+\frac{k\left(  k+11\right)  }{2}-\sum_{j=0}^{k-1}\left[  \frac{p+j}%
{2}\right]
\end{align*}
Now $\left[  \frac{8+j}{2}\right]  \leq\left[  \frac{7+k}{2}\right]  $ for
$j\leq k$, and $\left[  \frac{7+k}{2}\right]  \leq\frac{7+k}{2}$, with
equality when $k$ is odd. Then
\[
\frac{k\left(  k+11\right)  }{2}-\sum_{j=0}^{k-1}\left[  \frac{p+j}{2}\right]
\leq\frac{k\left(  k+11\right)  -k\left(  7+k\right)  }{2}=2k
\]
which shows that $f\left(  8+k\right)  \geq9+k.$
\end{proof}

With this comparison, we are in situation of giving the best possible lower
bound for the number of lines of a $\left(  p,q\right)  $-weight graph :

\begin{proposition}
Let $G$ be a $\left(  p,q\right)  $-weight graph of a nilpotent Lie algebra
$\frak{g}$. Then
\[
q\geq\left(
\begin{array}
[c]{c}%
p\\
2
\end{array}
\right)  -\sum_{j=1}^{\left[  \frac{p}{2}\right]  }\left(  p-2j\right)
\]
\end{proposition}

These results establish necessary conditions to be satisfied by an arbitrary graph to be a weight graph for some weight system $R\frak{g}\left(T\right)$ of a nilpotent Lie algebra $\frak{g}$ :

\begin{theorem}
Let $G=\overline{G\left(  R\frak{g}\left(  T\right)  \right)  }$ be the
$\left(  p,q\right)  $-weight graph of a nilpotent Lie algebra. Then

\begin{enumerate}
\item  for any $n\geq2$ the graph $G$ is isomorphic to its $n^{th}$ power
$G^{n}$

\item  for any pair of points $p_{1},p_{2}$ of $G$ the distance is $d\left(
p_{1},p_{2}\right)  \leq2$

\item $q\geq\left(
\begin{array}
[c]{c}%
p\\
2
\end{array}
\right)  -\sum_{j=1}^{\left[  \frac{p}{2}\right]  }\left(  p-2j\right)  .$
\end{enumerate}
\end{theorem}

\begin{example}
If $G=G\left(  R\frak{g}\left(  T\right)  \right)  $ has six points,
then the bound $q\geq\left(  p-1\right)  $ would tell us that the graph $G$
$\ $has at least $5$ lines. Now the bound of the theorem tells that $G$ has at
least $9$ lines. Eliminating the non connected graphs, the number of candidate
graphs has been reduced in $66$, which constitutes a considerable reduction of
the possibilities.
\end{example}

In particular, we can obtain, as a direct consequence of this, a special case
of a result of Lie-theoretic nature shown in [Fa] :

\begin{proposition}
For a fixed dimension $n\geq3$, there exists only a finite number of
equivalence classes of weight systems satysfying condition 1.
\end{proposition}

\end{document}